\documentclass[11pt]{amsart}

\usepackage[latin1]{inputenc}
\usepackage{amsmath,amsfonts,amssymb}
\usepackage{amscd}
\usepackage{graphicx}

\newcommand{\Z}{{\mathbb{Z}}}
\newcommand{\N}{{\mathbb{N}}}

\newcommand{\Ext}{\operatorname{Ext}}

\newcommand{\Hom}{\operatorname{Hom}}
\newcommand{\rad}{\operatorname{rad}}

\newcommand{\gldim}{\operatorname{gldim}}

\renewcommand{\d}{\displaystyle}

\newtheorem{theorem}{Theorem}[section]
\newtheorem{prop}[theorem]{Proposition}
\newtheorem{cor}[theorem]{Corollary}
\newtheorem{lemma}[theorem]{Lemma}
\newtheorem{example}[theorem]{Example}
\newtheorem{remark}[theorem]{Remark}

\begin{document}

\title{Weighted locally gentle quivers and Cartan matrices}

\author{Christine Bessenrodt \and Thorsten Holm}

\address{~~\newline
Institut f\"ur Mathematik, Universit\"at Hannover,
Welfengarten 1, D-30167 Hannover, Germany}
\email{bessen@math.uni-hannover.de}

\address{
~~\newline
Department of Pure Mathematics,
University of Leeds,
Leeds LS2 9JT,
United Kingdom}
\email{tholm@maths.leeds.ac.uk}


\bigskip

\begin{abstract}
We study the class of weighted locally gentle quivers.
This naturally extends the class of gentle quivers and
gentle algebras, which have been
intensively studied in the representation theory of
finite-dimensional algebras, to a wider class of potentially
infinite-dimensional algebras. Weights on the arrows of these
quivers lead to gradings on the corresponding algebras.
For the natural grading by path lengths, any locally
gentle algebra is a Koszul algebra.

Our main result is a
general combinatorial formula for the
determinant of the weighted
Cartan matrix of a weighted locally gentle quiver. This determinant
is invariant under graded derived equivalences of the corresponding
algebras.
We show that this weighted Cartan determinant is a
rational function which is completely determined by the
combinatorics of the quiver, more precisely by the
number and the weight of certain oriented cycles.
This leads to combinatorial invariants of the graded derived
categories of graded locally gentle algebras.

By specializing to certain parameters we
reobtain and extend in this way results from
\cite{BH} on Cartan determinants of gentle algebras.
\end{abstract}

\maketitle


\section{Introduction}

\subsection{}
In the representation theory of finite-dimensional algebras,
gentle algebras occur naturally in various contexts, especially
in connection with tilting and derived equivalences \cite{AH},
\cite{AS}, \cite{V}.
The definition
of these algebras is purely combinatorial in terms of quivers
with relations.

Combinatorially, as well as algebraically, it is very natural to
drop the condition on a gentle quiver
that the corresponding algebra should be finite-dimensional,
leading us in this paper to study so-called locally gentle quivers.
However, many of the advanced techniques from the representation
theory of finite-dimensional algebras do not carry over from
gentle algebras to locally gentle algebras.
Locally gentle algebras are also of interest because they are
Koszul algebras (cf.\ Section~\ref{sec:duality}).

Our first aim in this article is to give a general combinatorial
formula for the determinant of the weighted Cartan matrix of a locally
gentle quiver
(see~\ref{Cartan-locally} for the definition of these
matrices, whose entries are power series).

The main motivation
comes from the finite-dimensional situation: the unimodular
equivalence class, and hence the determinant, of
the Cartan matrix is an invariant of the derived module category
of a finite-dimensional algebra. For the (finite-dimensional)
gentle algebras, the Cartan
determinant was completely described in~\cite{TH}, and the
unimodular equivalence class in~\cite{BH}, leading to
new combinatorial derived invariants for gentle and skewed-gentle
algebras.

In the course of the proof of our formula for the determinant of
the Cartan matrix of a locally gentle quiver we are
naturally led to consider weights on the arrows of the quivers.
So we introduce in this paper the class of weighted locally gentle
quivers, and define their weighted Cartan matrices.

To any weighted locally gentle quiver $(Q,I)$ there is a
corresponding locally gentle algebra~$A=KQ/I$.
Algebraically, the weights on the quiver correspond to gradings
of the algebra. If all arrows are set to be in degree~1, then
the corresponding locally gentle algebras are Koszul (see~\ref{sec-Koszul} below).

Our main object of study, the determinant
of the weighted Cartan matrix is known to be invariant under
graded derived equivalences.
So our results on weighted Cartan determinants will give rise
to combinatorial invariants of graded derived categories of
locally gentle algebras.

Our main result gives an explicit formula for the determinant
of the Cartan matrix of a weighted locally gentle quiver.
The numerator and denominator of this rational function
are completely determined by the combinatorics of the quiver,
more precisely by the minimal oriented cycles with full relations
(for the
numerator) and the minimal oriented cycles with no relations
(for the denominator). See~\ref{intro-mainthm} for a precise
statement.

Evaluating our formula at special values, we
reobtain as special cases the results on Cartan determinants
from the previous papers~\cite{TH} and~\cite{BH}, as outlined
in~\ref{intro-appl} below.

\subsection{} \label{intro-mainthm}
We now briefly describe and state our main result. For definitions
and notations we refer to Section~\ref{sec-weighted} below.
For any $w$-weighted locally gentle quiver $\mathcal{Q}=(Q,I)$
denote by $\mathcal{ZC}(\mathcal{Q})$ the set of (minimal)
oriented cycles with full relations
and by $\mathcal{IC}(\mathcal{Q})$ the set of (minimal)
oriented cycles with no relations.

For any path $p=\alpha_1\alpha_2\ldots\alpha_t$
in~$Q$ (no matter whether zero or non-zero
in~$(Q,I)$) we let $w(p)=\prod_{i=1}^t w(\alpha_i)$ be its
weight and $l(p)=t$ its length.
\medskip

\noindent
{\bf Main Theorem.} {\it Let $\mathcal{Q}=(Q,I)$ be a
locally gentle quiver with the generic weight function~$w$,
and let $C_{\mathcal{Q}}^w(x)$
be its weighted Cartan matrix.

Then the determinant of this Cartan matrix is a rational function which
is given by the formula
$$\det C_{\mathcal{Q}}^w(x) =
\frac{\prod_{C\in\mathcal{ZC}(\mathcal{Q})} (1-(-1)^{l(C)}{w(C)})}
{\prod_{C\in\mathcal{IC}(\mathcal{Q})} (1-{w(C)})}\:.
$$
}

\subsection{} \label{intro-appl}
From our main result we can draw several immediate
conclusions.

Setting all arrow weights equal to an indeterminate~$q$, we get the following
result for locally gentle algebras which generalizes
Corollary~1 of~\cite{BH} to the infinite-dimensional
situation.
\medskip

\noindent
{\bf Corollary 1.} {\em Let $\mathcal{Q}=(Q,I)$ be a locally gentle quiver,
with corresponding algebra~$A=KQ/I$. Then the determinant
of the $q$-Cartan matrix is
$$\det C_{\mathcal{Q}}(q) = \det C_{A}(q) =
\frac{\prod_{C\in\mathcal{ZC}(\mathcal{Q})} (1-(-q)^{l(C)})}
{\prod_{C\in\mathcal{IC}(\mathcal{Q})} (1-q^{l(C)})}.
$$
}

Specializing further to finite-dimensional gentle algebras
and setting $q=1$ we obtain
the following explicit formula for the Cartan determinant which
has first been proven in~\cite{TH}. The Cartan determinant
is of particular importance since it is an invariant of
the derived module category.
For a gentle quiver $\mathcal{Q}=(Q,I)$ denote by $ec(\mathcal{Q})$
the number of oriented cycles of even length
in~$Q$ with full relations, and by $oc(\mathcal{Q})$
the corresponding number of such cycles of odd length.
\medskip

\noindent
{\bf Corollary 2.} {\em Let $\mathcal{Q}=(Q,I)$ be a gentle quiver,
with corresponding finite-dimensional algebra~$A=KQ/I$.
Then the determinant of the
Cartan matrix of~$\mathcal{Q}$ (and hence of
$A$) is
 $$\det C_{\mathcal{Q}} = \left\{
\begin{array}{ll}
0 & \mbox{if $ec(\mathcal{Q})>0$} \\
2^{oc(\mathcal{Q})} & \mbox{else}
\end{array} \right..
$$
}

\subsection{} The paper is organized as follows.
In Section~\ref{sec-weighted} we give the definition of weighted locally gentle
quivers, and of their weighted Cartan
matrices. In Section \ref{sec:duality} we prove a duality result.
In the unweighted case, this can be seen as a special case of
Koszul duality; however, we provide a general elementary combinatorial
proof.
The core part of the paper is
Section~\ref{sec:mainproof} which contains the proof of the
main theorem. In subsection~\ref{finite} we first give a proof
of the main theorem on weighted Cartan determinants
for the case where the weighted locally gentle quiver does not contain
a cycle with no relations. This is basically the finite-dimensional
situation considered in \cite{BH}, but now generalized to weighted
quivers. Actually, we get a more precise result about
unimodular normal forms for the weighted Cartan matrices.
Subsection \ref{reduce-ZC} contains a crucial reduction
step; we show how, under certain conditions, one can reduce
the number of cycles with full relations in the quiver,
and at the same time precisely control the transformations on the
weighted Cartan matrices. Using this reduction result, we will then
in subsection~\ref{dual} prove the main theorem
for `most' weighted locally gentle quivers. Actually, by
the reduction via~\ref{reduce-ZC} we inductively get a quiver with
no cycles with full relations. Hence its dual weighted locally
gentle quiver has no cycles with no relations, i.e., we are
in the finite-dimensional situation of~\ref{finite}. Applying
the duality result~\ref{duality-thm} then finishes the proof
of the main theorem.
In subsection~\ref{explicit-example} we illustrate the above
arguments by going through an explicit example.

However, there are very special quivers for which these
arguments do not work. These so-called critical locally
gentle quivers are introduced and dealt with in subsection~\ref{critical},
where we compute their weighted Cartan determinant.
The critical locally gentle quivers have an interesting combinatorial
interpretation in terms of certain configurations of $2n$-polygons
(where $n$ is the number of vertices of the quiver). We explain
this connection in detail in subsection \ref{configurations}.
In particular, as a consequence of the Harer-Zagier formula
\cite{Harer-Zagier},
critical locally gentle quivers exist only with an even number of
vertices.

\subsection{Acknowledgement} We are grateful to several
people for very helpful discussions about the topics of this paper.
In particular, we thank Joe Chuang, Peter J{\o}rgensen
and Rapha\"el Rouquier for patiently answering questions about Koszul
algebras and graded derived equivalences, and Richard Stanley for
pointing out to us the
references \cite{Goulden-Nica}, \cite{Harer-Zagier}.

\section{Weighted locally gentle quivers} \label{sec-weighted}

In this section we will eventually introduce the weighted locally
gentle quivers occurring in the title. Before, we give a quick
review of the basics on quivers, Cartan matrices and gentle
quivers in the ordinary (non-weighted and finite-dimensional)
setting.

\subsection{Quivers with relations}

Algebras can be defined naturally from a combinatorial setting
by using directed graphs. A {\em quiver} is a directed graph
with finitely many vertices and arrows.

For any field~$K$, we can define the
path algebra~$KQ$. It has as
basis the set of all oriented paths in~$Q$. The multiplication
is defined
by concatenation of paths. More precisely, for a path~$p$
in~$Q$ let $s(p)$ denote its start vertex and $t(p)$ its
end vertex. The product in~$KQ$ of two paths $p$ and~$q$ is defined
to be the concatenated path~$pq$ if $t(p)=s(q)$, and zero
otherwise. Note that our convention is to write paths from
left to right.

Such a path algebra~$KQ$ is finite dimensional precisely when
$Q$ does not contain an oriented cycle.

More general algebras can be obtained by introducing relations on a
path algebra.
An ideal $I\subseteq KQ$ is called admissible if $I\subseteq
\rad^2(KQ)$ where $\rad(KQ)$ is the radical of the algebra~$KQ$.

A famous theorem of P.\ Gabriel
states that if~$K$ is algebraically closed, any finite-dimensional
$K$-algebra is Morita equivalent to a factor algebra
$KQ/I$ where $I$ is an admissible ideal.

So for most contexts within the representation theory of
finite-dimensional algebras
it suffices to consider algebras of the form
$KQ/I$, often referred to as {\em quivers with relations}.

\subsection{Ordinary Cartan matrices}
Let $A=KQ/I$ be finite-dimensional, where $K$ is any field.
By a slight abuse of notation we
identify paths in the quiver~$Q$ with their
cosets in~$A$. Let $Q_0$ denote the set of vertices of~$Q$.
For any $i\in Q_0$ there exists a path $e_i$ of length zero.
These are primitive orthogonal idempotents in~$A$,
the sum $\sum_{i\in Q_0} e_i$ is the unit element in~$A$.
In particular we get $A=1\cdot A=\oplus_{i\in Q_0} e_iA$,
hence the (right) $A$-modules
$P_i:=e_iA$ are the indecomposable projective
$A$-modules.

The (ordinary) {\em Cartan matrix} $C=(c_{ij})$ of a
finite-dimensional algebra $A=KQ/I$
is the $|Q_0|\times |Q_0|$-matrix defined by setting
$c_{ij}:=\dim_K \Hom_A(P_j,P_i)$.
Any homomorphism $\varphi:e_jA\to e_iA$ of right $A$-modules
is uniquely determined by $\varphi(e_j)\in e_iAe_j$,
the $K$-vector space generated by all paths in~$Q$ from vertex~$i$
to vertex~$j$, which are nonzero in~$A=KQ/I$.
In particular, we have $c_{ij}=\dim_Ke_iAe_j$. In this way,
computing entries of the Cartan matrix for $A=KQ/I$ is the same
as counting
paths in the quiver~$Q$ which are nonzero in~$A$.

This is the key viewpoint in this paper, enabling us to
obtain results on the representation-theoretic Cartan invariants
by purely combinatorial methods.

\subsection{Gentle quivers}

A pair $(Q,I)$ consisting of a quiver~$Q$ and an admissible
ideal~$I$ in the path algebra~$KQ$ is called
a {\em special biserial quiver} if it satisfies the
following axioms (G1)-(G3).
\smallskip

(G1) The corresponding algebra $A=KQ/I$ is finite-dimensional.
\smallskip

(G2) Each vertex of~$Q$ is starting point of at most two arrows,
and end point of at most two arrows.
\smallskip

(G3) For each arrow $\alpha$ in~$Q$ there is at most one arrow
$\beta$ such that $\alpha\beta\not\in I$, and at most one arrow
$\gamma$ such that $\gamma\alpha\not\in I$.
\smallskip

Gentle quivers form a subclass of the class of
special biserial quivers.
\smallskip

The pair $(Q,I)$ as above is called {\em gentle} if it is special
biserial, i.e. (G1)-(G3) hold,
and in addition the following axioms hold.
\smallskip

(G4) The ideal $I$ is generated by paths of length~2.
\smallskip

(G5) For each arrow $\alpha$ in~$Q$ there is at most one arrow
$\beta'$ with $t(\alpha)=s(\beta')$ such that $\alpha\beta'\in I$,
and there is at most one arrow $\gamma'$ with
$t(\gamma')=s(\alpha)$ such that $\gamma'\alpha\in I$.
\smallskip

An algebra~$A$ over the field~$K$
is called {\em gentle}
if $A$ is Morita equivalent to an algebra~$KQ/I$ where
$(Q,I)$ is a gentle quiver.

\subsection{Locally gentle quivers} \label{locally}

Gentle algebras have been introduced and studied intensively in
the representation theory of finite-dimensional algebras. This explains
the occurrence of axiom~(G1) on finite-dimensionality above.
However, from a
combinatorics point of view this axiom does not seem to be
natural; in fact any 'combinatorial' statement or proof
about gentle algebras
(like e.g. the ones about Cartan determinants and normal forms
in \cite{TH}, \cite{BH}) can be hoped to have a counterpart in the
absence of~(G1).

In this paper we will be concerned also with infinite-dimensional
algebras arising from these quivers with relations. This leads us to
the following definition.

A {\em locally gentle quiver} is a
pair $(Q,I)$ consisting of a quiver~$Q$ and an admissible
ideal~$I$ in the path algebra~$KQ$ satisfying
(G2)-(G5) above.

To any such locally gentle quiver $(Q,I)$ there is attached
a {\em locally gentle algebra}~$A=KQ/I$.
\smallskip

Simple examples of locally gentle algebras
are given by the polynomial ring~$K[X]$ in one indeterminate
over a field~$K$, or the (non-commutative) algebra
$K\langle X,Y\rangle/(X^2,Y^2)$.

\subsection{Cartan matrices for locally gentle quivers}
\label{Cartan-locally}

Since a locally gentle algebra need not be finite-dimensional,
the usual definition of a Cartan matrix (where the entries are obtained
by just counting the number
of non-zero paths) no longer makes sense.

Instead, as in~\cite{BH}, one can look at
a refined version of the Cartan matrix, where
instead of just counting paths in~$(Q,I)$, we now count
each path (which is non-zero in the algebra~$A=KQ/I$) of length~$n$
by~$q^n$, where $q$ is an indeterminate.

More precisely,
the path algebra~$KQ$ is a graded algebra, with grading given by
path lengths. Since $I$ is homogeneous, the factor algebra~$A=KQ/I$
inherits this grading. So
the morphism spaces $\Hom_A(P_j,P_i)\cong e_iAe_j$ become
graded vector spaces.
For any vertices~$i$ and~$j$ in~$Q$
let $e_iAe_j= \oplus_{n} (e_iAe_j)_n$ be the
graded components.

Now let $q$ be an indeterminate.
For a locally gentle quiver ${\mathcal Q}=(Q,I)$
the $q$-Cartan matrix $C_{{\mathcal Q}}=(c_{ij}(q))$ of~${\mathcal Q}$
is defined as the matrix with entries
$c_{ij}(q):= \sum_n \dim_K (e_iAe_j)_n\, q^n$
in the ring of power series $\mathbb{Z}[[q]]$.
\medskip

In other words, the entries of the $q$-Cartan matrix are the
Poincar\'{e} polynomials of the graded homomorphism spaces
between projective modules of the corresponding algebra
$A=KQ/I$.

Clearly, specializing $q=1$ gives back the usual Cartan matrix
(i.e., we forget the grading).

\subsection{Weighted locally gentle quivers}
\label{weight-def}

We now introduce the most general class of quivers to be studied in
this paper.
\smallskip

A {\em $w$-weighted locally gentle quiver}
is a locally gentle quiver $(Q,I)$ together
with a {\em weight function} $w: Q_1\to R$
on the arrows of~$Q$ into a ring~$R$ with~1.

The weight function is extended to all paths in~$Q$
by setting $w(p)=1$ for the trivial paths of length 0, and
$w(p)=\prod_{i=1}^t w(\alpha_i)$
for a  path $p=\alpha_1\cdots \alpha_t$ with
$\alpha_1, \ldots, \alpha_t \in Q_1$.
In later situations, the weight function will be
further restricted so that corresponding weighted counts make sense.

\smallskip

Note that the special case of choosing the weight function
$w: Q_1 \to \Z[q]$, $\alpha \mapsto q$ for all arrows $\alpha$,
induces the weight by length  on the
paths in the quiver. This has been used before for the $q$-Cartan matrix
of the ordinary locally gentle quivers, as defined in~\ref{locally}.
We will also refer to these as {\em $q$-weighted}
locally gentle quivers.

\smallskip

The main case for us is the {\em generic weight function}
$w:Q_1 \to \Z[x_e \mid e\in Q_1]$ into the polynomial ring
$\Z[x_e \mid e\in Q_1]$ given by mapping each arrow $e\in Q_1$
to the corresponding indeterminate~$x_e$.
Of course, specializing all $x_e$ to $q$ leads to the previous
case.
If $t$ is a further indeterminate,
substituting $x_et$ for $x_e$, for all $e\in Q_1$,
gives a weight function on the paths that
shows explicitly both the generic weight and the weight by length (here via
$t^{l(p)}$ for a path~$p$).
For later purposes, we may as well substitute other monomials from a
polynomial ring $\Z[y_1, \ldots,y_k]$ for the~$x_e$'s.
For example, the weights may be of the form $q^{m_e}t$, $m_e\in \N$,
for $e\in Q_1$; with this choice we may keep track at the same time
of the length and a further integer weight of a path.

\subsection{The Cartan matrix of a weighted locally gentle quiver}
\label{Cartan-weight}

Let ${\mathcal Q}=(Q,I)$ be a locally gentle quiver with
the generic weight function $w:Q_1\to \Z[x_e\mid e\in Q_1]$
as above.
We define a weighted Cartan matrix $C_{\mathcal{Q}}^w(x)$
(where $x$ stands for $(x_e)_{e\in Q_1}$)
as for locally gentle quivers in \ref{Cartan-locally}
by counting non-zero paths according to their weights,
where here instead of the lengths we take
into account the weights on the arrows.
Thus, a non-zero path $p$ in the corresponding algebra
$A=KQ/I$ gives a contribution~$w(p)$.
The corresponding Cartan matrix is then defined
over the ring of power series $\Z[[x_e\mid e\in Q_1]]$.
For any vertices $i,j \in Q_0$ the
corresponding entry in $C_{\mathcal{Q}}^w(x)$ is set to be
$$c_{ij}(x):=\sum_{p} {w(p)}$$
where the sum is taken
over all non-zero paths $p$ in $(Q,I)$ from $i$ to~$j$.
\smallskip

Note that if in the weight function all variables are
specialized to~$q$ on~$Q_1$, i.e.,
if it induces the weighting by path lengths,
then we reobtain the $q$-Cartan matrix $C_{\mathcal{Q}}(q)$
of the quiver.


\section{Duality} \label{sec:duality}

In this section we will define to any weighted locally gentle
quiver its dual and prove a fundamental duality result for
their Cartan matrices which will be crucial later in the proof
of the main theorem.

\subsection{The dual weighted locally gentle quiver}
\label{subsec:dual}
Let $\mathcal{Q}=(Q,I)$ be a $w$-weighted locally gentle
quiver. The {\em dual} $w$-weighted locally gentle quiver
is defined as $\mathcal{Q}^{\#}:=(Q,I^{\#})$
where the relation set $I^{\#}$ is characterized by the
following property. For any
arrows $\alpha,\beta$ in~$Q$ with $t(\alpha)=s(\beta)$
we have $\alpha\beta\not\in I$ if and only if
$\alpha\beta\in I^{\#}$.

Note that the dual $\mathcal{Q}^{\#}$ has the same underlying
quiver~$Q$ and the same weight function~$w$ as~$\mathcal{Q}$.

We would like to point out that the above duality
does not preserve finite-dimensionality for the corresponding algebras,
as the following easy example shows.
\smallskip

Take $Q$ to be the quiver with one vertex and one loop~$\epsilon$.
Let $I$ be generated by~$\epsilon^2$. Then $(Q,I)$
is gentle, giving a
finite-dimensional gentle algebra~$A=KQ/I$. For the dual
quiver $(Q,I^{\#})$ we have $I^{\#}=0$, so the corresponding
algebra $A^{\#}=KQ/I^{\#}$ is infinite-dimensional (more precisely,
it is a polynomial ring in one generator).
\smallskip

This explains that for being able to use duality we really
have to broaden our perspective and study infinite-dimensional
situations as well.

\subsection{Duality on Cartan matrices} \label{sec-duality}

The aim is to prove the following fundamental result on
the Cartan matrices of weighted locally gentle quivers.

\begin{prop} \label{duality-thm}
Let $\mathcal{Q}=(Q,I)$ be a locally
gentle quiver with the generic weight function~$w$,
and let $\mathcal{Q}^{\#}=(Q,I^{\#})$ be its dual.
Then the following holds for their weighted Cartan matrices:
$$C_{\mathcal{Q}}^w(x)\cdot  C_{\mathcal{Q}^{\#}}^w(-x) = E_{|Q_0|}$$
where $E_{|Q_0|}$ denotes the identity matrix.
\end{prop}

\proof For $i,j\in Q_0$, we consider the $(i,j)$-entry of the product
$C_{\mathcal{Q}}^{w}(x)\cdot  C_{\mathcal{Q}^{\#}}^{w}(-x)$.
Any contribution to this comes from a path
$\hat p=pp'$ from $i$ to~$j$ in~$KQ$ such that the path~$p$ is non-zero
in~$A=KQ/I$ and the path~$p'$ is non-zero in~$A^{\#}=KQ/I^{\#}$.
Let $p=p_0p_1\ldots p_l$, $p'=p_0'p_1'\ldots p_{l'}'$ with
$p_r,p_r' \in Q_1$,
except that possibly $p_0=e_i$ if $l=0$ or $p_0'=e_j$ if $l'=0$.
The contribution coming from this path is then
$w(p)w(p')(-1)^{l(p')}$.
By definition of the dual quiver,
$p_lp_0'\not\in I$ or  $p_lp_0'\not\in I^{\#}$, and we can have both only
if (at least) one of the paths is trivial.
If $p_lp_0'\not\in I$ and $l'>0$, set $p''=p_1'\ldots p_{l'}'$;
then the same path $\hat p$ in $KQ$ also gives
a contribution coming from its factorization $\hat p=(pp_0')p''$,
which is $w(p)w(p')(-1)^{l(p')-1}$ and thus cancels
with the previous contribution. Similarly, we obtain a
cancelling contribution if $p_lp_0'\not\in I^{\#}$ and $l>0$
by shifting the factorization one
place to the left rather than to the right.
Note that for any fixed path $\hat p$ in $KQ$ of positive length
which gives a contribution
we have exactly two factorizations as a product of a non-zero path
in~$A$ and a non-zero path in~$A^{\#}$, as described above.
As the corresponding contributions cancel, it only remains to consider
the trivial paths at each vertex; these give a contribution~$1$,
and thus the matrix product is the identity matrix, as claimed.
\qed

\medskip

\begin{remark}\label{rem-dual}
{\rm
This proposition implies immediately that if the main theorem
holds for a weighted locally gentle quiver ${\mathcal Q}$ then it
is also true for its dual.
To see this, we first note that clearly cycles with full relations
and cycles without relations are interchanged when we dualize
the quiver.
Furthermore, evaluating the Cartan matrix at $-x$ rather than at~$x$
means that the weight of a cycle $C$ of length~$l(C)$ is changed
from $w(C)$ to $(-1)^{l(C)}w(C)$. Hence, if the determinant formula
holds for the Cartan matrix of~${\mathcal Q}$, then the duality formula above gives
$$
\begin{array}{rcl}
\det C_{\mathcal{Q}^{\#}}^w(x) &=&
\d (\det C_{\mathcal{Q}}^w(-x))^{-1}\\[10pt]
& =&
\d \frac{\prod_{C\in\mathcal{IC}(\mathcal{Q})} (1-(-1)^{l(C)}{w(C)})}
{\prod_{C\in\mathcal{ZC}(\mathcal{Q})} (1-{w(C)})}\\[10pt]
&=&
\d \frac{\prod_{C\in\mathcal{ZC}(\mathcal{Q}^{\#})} (1-(-1)^{l(C)}{w(C)})}
{\prod_{C\in\mathcal{IC}(\mathcal{Q}^{\#})} (1-{w(C)})}\,.
\end{array}
$$
In particular, we observe that if the main theorem holds for all
weighted locally gentle quivers which do not have cycles with no relations (i.e.,
those corresponding to finite-dimensional algebras), then it
also holds for all weighted locally gentle quivers which do not have cycles with full
relations.
}
\end{remark}

\begin{example} \label{ex-duality}
{\rm We consider the following weighted locally gentle
quiver $\mathcal{Q}$ (where relations are indicated by
dotted lines)


\begin{center}
\includegraphics[scale=0.5]{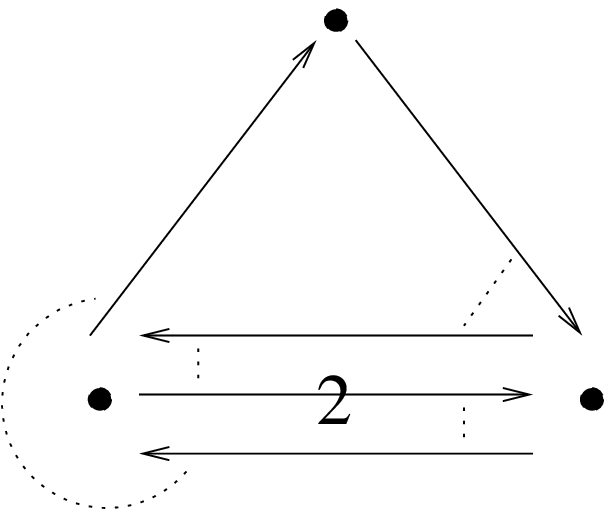}
\end{center}


A number $m$ attached to an arrow denotes the weight $q^mt$, where we
omit attaching 1's for the weight~$qt$; again, $q,t$ are indeterminates.
The vertices are denoted $1,2,3$ in a clockwise order,
with $1$ being the vertex at the top.
Then the weighted Cartan matrix has the following form
$$C_{\mathcal{Q}}(q,t)=
\left(\begin{array}{ccc}
\frac{1}{1-q^6t^5} & \frac{qt+q^4t^3}{1-q^6t^5} &
\frac{q^2t^2+q^5t^4}{1-q^6t^5} \\
\frac{q^2t^2+q^5t^4}{1-q^6t^5} & \frac{1+q^3t^2+q^3t^3+q^6t^5}{1-q^6t^5} &
\frac{2qt+q^4t^3+q^4t^4}{1-q^6t^5} \\
\frac{qt+q^4t^3}{1-q^6t^5} & \frac{q^2t+q^2t^2+2q^5t^4}{1-q^6t^5} &
\frac{1+q^3t^2+q^3t^3+q^6t^5}{1-q^6t^5}
\end{array}
\right)
$$
The dual weighted locally gentle quiver $\mathcal{Q}^{\#}$ has
the form


\begin{center}
\includegraphics[scale=0.5]{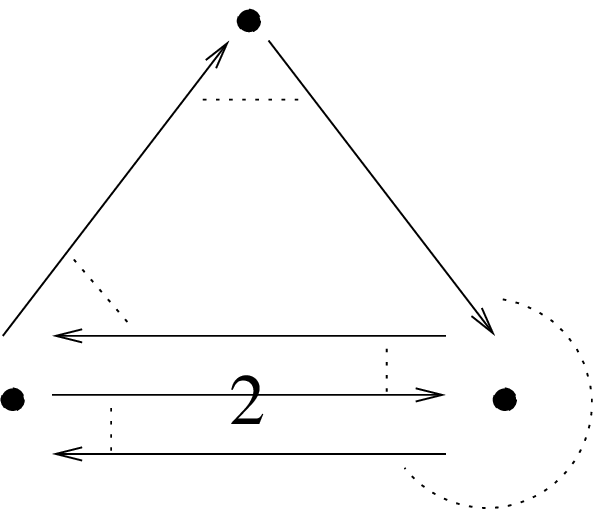}
\end{center}

Note that this dual quiver now has one minimal oriented cycle
with full relations of length 6, and no cycles with no relations.
The weighted Cartan matrix of $\mathcal{Q}^{\#}$ has the following
form
$$C_{\mathcal{Q}^{\#}}(q,t)=
\left(\begin{array}{ccc}
1+q^6t^5 & qt+q^4t^3 & q^2t^2+q^5t^4 \\
q^2t^2+q^5t^4 & 1+q^3t^2 & 2qt+q^4t^3 \\
qt+q^4t^3 & q^2t & 1+q^3t^2
\end{array}
\right)
$$
In this example, we have specialized the generic weights
$x_e$ to $q^{m_e}t$,
where $m_e$ is the weight number attached to the arrow~$e$ in the figure.
Thus, the change from $x$ to $-x$ for the Cartan matrix of the
dual quiver in the duality result corresponds to a change $t$ to $-t$
in this situation.
We leave to the reader the straightforward (though tedious)
verification that indeed
the product $C_{\mathcal{Q}}(q,t)\cdot C_{\mathcal{Q}^{\#}}(q,-t)$
equals the identity matrix, as predicted by Proposition~\ref{duality-thm}.
}
\end{example}

\subsection{Locally gentle algebras are Koszul} \label{sec-Koszul}

The above duality result Proposition~\ref{duality-thm}
is reminiscent of similar formulae for Koszul algebras,
so-called {\em numerical Koszulity criterion}, see~\cite[2.11.1]{BGS}.
In this subsection we briefly clarify this connection
by indicating that the locally gentle algebras,
with the natural grading by path lengths (i.e., all arrows are
of degree~1),
are actually Koszul algebras. However, if at least one of the
weights is~$>1$ then the weighted locally gentle algebra
can not be Koszul since any Koszul algebra has to be generated
by its degree~1 component (see~\cite[2.3.1]{BGS}).
In particular, our duality result above (which holds
for arbitrary weights) is not just a special case of
Koszul duality.

For background and more details on Koszul algebras we refer
to~\cite{BGS}, especially Section~2.

Recall that a positively graded algebra $A=\oplus_{i\ge 0} A_i$
is Koszul if $A_0$ is semisimple and if the graded left
$A$-module $A_0$ admits a graded projective resolution
$\ldots \to P^2\to P^1\to P^0\to A_0\to 0$ such that each $P^i$
is generated by its component in degree~$i$, i.e., $P^i=AP_i^i$.

Also recall that on any graded $A$-module $M=\oplus_{i\ge 0} M_i$
there are shifts defined by $(M\langle n\rangle)_i=M_{i-n}$.

The following observation completely describes the graded
projective resolutions of simple modules for locally gentle algebras.
We leave the details of the straightforward verification to the reader.

\begin{prop} \label{prop:resolution}
Let $(Q,I)$ be a locally gentle quiver, with
corresponding algebra $A=KQ/I$. Consider $A$ with the natural
grading given by path lengths.
For any vertex $i$ in~$Q$, consider the (at most) two
paths in~$Q$ with full relations starting in~$i$;
denote the vertices on these paths by $i,i_1,i_2,\ldots$
and by $i,j_1,j_2,\ldots$, respectively.
Then the corresponding simple module
$E_i$ has a graded projective resolution of the form
$$  \ldots\to (P^{i_2}\oplus P^{j_2})\langle -2\rangle
 \to (P^{i_1}\oplus P^{j_1})\langle -1\rangle
\to P^i \to E_i\to 0.
$$
In particular, with the grading by path lengths, any
locally gentle algebra $A$ is Koszul.
\end{prop}

As a consequence of the above proposition we get that
the Koszul dual $E(A)=\Ext_{A}(A_0,A_0)$ is canonically
isomorphic to the opposite algebra of the quadratic dual
$A^{!}$ (\cite[Theorem 2.10.1]{BGS}). But for a locally
gentle quiver ${\mathcal Q}=(Q,I)$ with algebra $A=KQ/I$,
the latter opposite quadratic dual is just the algebra
$A^{\#}:=KQ/I^{\#}$ given by the dual
locally gentle quiver ${\mathcal Q}^{\#}=(Q,I^{\#})$
defined in Section~\ref{subsec:dual}, having the 'opposite'
relations. (For a general definition of the quadratic dual, see
\cite[ Definition 2.8.1]{BGS}.)

Then our duality result Proposition~\ref{duality-thm},
for the case of the grading given by path lengths, is a special case
of the numerical Koszulity criterion \cite[Lemma 2.11.1]{BGS}.
\smallskip

A direct consequence of Proposition~\ref{prop:resolution}
is the following homological property of
locally gentle algebras. As usual, $\gldim(A)$ denotes the global
dimension of an algebra. Recall that by $\mathcal{ZC(Q)}$ we denote
the set of minimal oriented cycles with full relations.

\begin{cor} Let ${\mathcal Q}=(Q,I)$ be a locally gentle quiver,
with corresponding algebra $A=KQ/I$. Then
$$\gldim(A) < \infty ~~~ \Longleftrightarrow ~~~
\mathcal{ZC(Q)} = \emptyset.$$
\end{cor}


\section{Proof of the main result} \label{sec:mainproof}

In this section we will give a complete proof of our main result
which we recall here for the convenience of the reader.
\medskip

\noindent
{\bf Main Theorem.}
{\it Let $\mathcal{Q}=(Q,I)$ be a
locally gentle quiver with the generic weight function~$w$,
and let $C_{\mathcal{Q}}^w(x)$
be its weighted Cartan matrix.

Then the determinant of this Cartan matrix of $\mathcal{Q}$ is
given by the formula
$$\det C_{\mathcal{Q}}^w(x) =
\frac{\prod_{C\in\mathcal{ZC}(\mathcal{Q})} (1-(-1)^{l(C)} w(C) )}
{\prod_{C\in\mathcal{IC}(\mathcal{Q})} (1-{w(C)})}.
$$
}
\medskip

The proof will consist of three main steps.
First, we give a proof for the finite-dimensional case.
Secondly, a reduction to the case where the quiver has no
oriented cycles with full relations. Using Remark~\ref{rem-dual}
this inductively proves the main theorem for 'most' quivers.
Finally, we have to deal with certain so-called critical quivers
separately.


\subsection{The finite-dimensional case} \label{finite}
We first show how to prove the main result in the special case where
a weighted locally gentle quiver $(Q,I)$ has no non-zero infinite
paths.
Note that  this corresponds to
the algebra $A=KQ/I$ being finite-dimensional, i.e.,
we are in the situation of a weighted gentle quiver.
The $q$-weighted special case of the following result has already
been proven in \cite{BH}.

\begin{prop} \label{prop-finite}
Let $\mathcal{Q}=(Q,I)$ be a gentle quiver,
and let $w$ be the generic weight function on the quiver.
Then the weighted Cartan matrix $C_{\mathcal{Q}}^w(x)$
can be transformed by unimodular
elementary operations over $\Z[x]=\Z[x_1,\ldots,x_{|Q_0|}]$
into a diagonal matrix with entries $1-(-1)^{l(C)}w(C)$,
for each $C\in \mathcal{ZC}(\mathcal{Q})$, and all further
diagonal entries being~1.
\\
In particular, for the
determinant of
the weighted Cartan matrix we have
$$\det C_{\mathcal{Q}}^w(x) = \prod_{C\in\mathcal{ZC}(\mathcal{Q})}
(1-(-1)^{l(C)} w(C)).$$
\end{prop}

\proof The proof is analogous to the proof for the $q$-Cartan
matrix of (finite-dimensional) gentle algebras given in
\cite{BH}. However, we have to take the weights into account,
so that we should include a proof here (although we shall be
brief at times; for details we refer to \cite{BH}).

Since $\mathcal{Q}$ is gentle,
we can do a similar reduction as in \cite[Lemma 3.1]{BH}
and as at the beginning of the proof of \cite[Theorem 3.2]{BH},
now adapted to the weighted case.
After the reduction, we can assume that $Q$ contains a vertex $v=v_1$
of degree~2,
with incoming arrow {$p_0:v_0\to v_1$} and outgoing arrow $p_1:v_1\to v_2$
with $p_0p_1=0$ in $A$, where $v_0\neq v_1\neq v_2$. \\
As the quiver is gentle, there exists a unique
path~$p$ starting with $p_0$ such that the product of any two
consecutive arrows on~$p$ is in~$I$.
As we have already gone through the reduction steps,
this path is an oriented cycle $\mathcal C$ with full relations
returning to~$v_0$.
In this case we set $p$ to be just one
walk around the cycle. So we have defined  a finite path
$p=p_0p_1\ldots p_s$ (with full relations). Let
$v_0,v_1,\ldots,v_s,v_{s+1}=v_0$ denote the vertices on the path~$p$.

We shall perform elementary row operations on the Cartan matrix
$C^w_{\mathcal{Q}}(x)$. Denote the row corresponding to a vertex
$u$ of~$Q$ by~$z_u$.
Then consider the following linear combination of rows
$$Z:=z_{v_1} - {w(p_1)}z_{v_2} + {w(p_1)w(p_2)}z_{v_3}
-+\ldots \hskip4cm$$
$$\hskip1cm \ldots+ (-1)^{s-1}{w(p_1)\cdots w(p_{s-1})z_{v_{s}}}
+(-1)^{s}{w(p_1)\cdots w(p_{s})z_{v_{s+1}}}.
$$
We replace the row $z_{v_1}$ by $Z$ and get a new matrix
$\widetilde{C}$.
The crucial observation is that in the alternating sum~$Z$
many parts cancel. In fact, at any vertex~$v_i$ on~$p$
there is a bijection between the non-zero paths starting in~$v_i$
but not with~$p_i$, and the non-zero paths starting with~$p_{i-1}$.
The `scalar' factors in~$Z$ are just chosen appropriately
so that the corresponding contributions in the weighted Cartan
matrix cancel.
Hence, in $Z$ only those contributions could survive
coming from paths starting in~$v_1$, but not with~$p_1$.
But by the choice of $v_1$ there are no such paths except the trivial one.
For all other rows in $C^w_{\mathcal{Q}}(x)$, note that no path
from some vertex $\tilde{v}\neq v_1$ can involve~$p_1$.
(In fact, there is only one incoming arrow in~$v_1$.)

As $p=\mathcal C$ is an oriented cycle with full relations,
we also have a final trivial bijection
$\{p_0\}\to \{e_{v_1}\}$ where the weight is reduced by $w(p_0)$.
By the remark above, almost everything cancels in the new row~$Z$
apart from the one term in the column to $v_1$ which comes from
the trivial path at $v_1$ and its multiplication by all the weights, i.e.
$$1+(-1)^{s} {w(p_1)\cdots w(p_{s})} {w(p_0)} =
1-(-1)^{s+1} {w(p_1 \ldots p_{s}p_0)}=1-(-1)^{l({\mathcal C})}
w({\mathcal C})\:.$$

As in the situation with the length weight function in~\cite{BH},
we then use the corresponding operation on the columns labelled by the vertices
on the cycle~$\mathcal C$, but now in counter-clockwise order, i.e., we set
$v_{s+2}=v=v_1$ and replace the column $s_{v_1}$ by
$$S:=s_{v_1} - {w(p_0)}s_{v_0} + {w(p_s)w(p_0)}s_{v_s}
-+\ldots \hskip4cm$$
$$\hskip1cm \ldots+ (-1)^{s-1}w(p_{3})\cdots w(p_s)w(p_0)s_{v_{3}}
+(-1)^{s}{w(p_{2})\cdots w(p_s)w(p_0)s_{v_{2}}}.
$$
Again, this amounts to the desired cancellation of terms. Thus,
by also ordering vertices so that
$v$ corresponds to the first row and column of the Cartan matrix,
we have altogether transformed $C^w_{\mathcal Q}(x)$ to a matrix
of the form
$$\left(
\begin{array}{cccc}
1-(-1)^{l({\mathcal C})}w({\mathcal C}) & 0 & \cdots & 0 \\
0 \\
\vdots & & C' \\
0\\
\end{array}
\right)
$$
Here $C'$ is the  weighted Cartan matrix
for the gentle quiver $\mathcal Q'$ obtained
from ${\mathcal Q}$ by removing $v$ and  the arrows incident with $v$
(and removing the corresponding relations)
and restricting the weight function to $Q_1'$ (so this
is the generic weight function $w'$ for $\mathcal Q'$).
Note that in comparison with~$\mathcal Q$, the quiver $\mathcal Q'$ has
one vertex less and one cycle with full relations
less (namely $\mathcal C$).
Now by induction, the result holds for
$C'=C^{w'}_{{\mathcal Q'}}(x')$,
and hence the result for $C^w_{\mathcal Q}(x)$ follows immediately.
\qed

\subsection{Reducing the zero cycles} \label{reduce-ZC}

The aim of this section is to prove a technical
result which is actually the main reduction step for the proof
of the main theorem. It describes a combinatorial
procedure for reducing the number of oriented cycles with full
relations, leading to a new weighted locally gentle quiver.
The crucial aspect is that we can control
the transformations on the determinants of
the weighted Cartan matrices in this process.
However, in the reduction procedure we are going to replace
two consecutive arrows by one new arrow, with weight equal to the
product of the weights of the former arrows. But this process
changes the lengths of cycles, so that in the following
result the weighted determinant of the new quiver
as a function in~$x$ is not quite with respect to
the generic weight function on the new quiver.

\begin{prop} \label{prop-reduce}
Let $\mathcal{Q}=(Q,I)$ be a locally gentle quiver
which contains a minimal oriented cycle~$C$ with full relations;
let $w$ be the generic weight function.
Assume there exists a vertex $v_1$ on~$C$ such that only two arrows
of~$C$ are incident with~$v_1$.
Let $p_1$ be the arrow on~$C$
with starting point~$v_1$.
Assume that there exists an
incoming arrow  $q_1$ at~$v_1$ which does not belong to~$C$.
We define a new $\bar{w}$-weighted locally gentle quiver
$\bar{\mathcal{Q}}=(\bar{Q},\bar{I})$ as follows.
The vertices are the same
as in~$Q$, the arrows $p_1$ and~$q_1$ are removed, and replaced
by one arrow~$\bar{p}$ with $s(\bar{p})=s(q_1)$, $t(\bar{p})=t(p_1)$.
The weight function $\bar w$ is set to
be $\bar w(\bar{p}):=w(q_1)w(p_1)$ on the new arrow,
all other weights are the same as for~$w$.

If $C$ is the only minimal oriented cycle with full relations
attached to~$v_1$ then
$$\det C_{\mathcal{Q}}^w(x) = (1-(-1)^{l(C)} w(C) )\cdot
\det C_{\bar{\mathcal{Q}}}^{\bar w}(x).$$

If we have a second minimal oriented cycle~$C'$ with full relations
attached to~$v_1$, then
$$\det C_{\mathcal{Q}}^w(x) =
(1-(-1)^{l(C)} w(C))(1-(-1)^{l(C')} {w(C')})\cdot
\det C_{\bar{\mathcal{Q}}}^{\bar w}(x).$$
\end{prop}

The following figure illustrates the situation and the statement
of the above proposition.

\bigskip


\begin{center}
\includegraphics[scale=0.5]{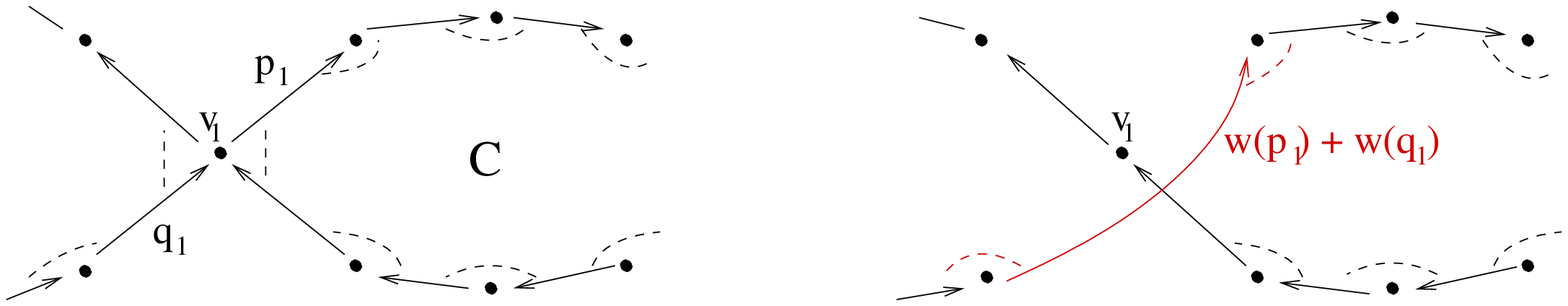}
\end{center}


\bigskip

\proof Analogous to the proof of Proposition~\ref{prop-finite}
we perform row operations along the cycle~$C$. Let
$C=p_1p_2\ldots p_k$ be the arrows, and
$v_1,v_2,\ldots,v_k$ the vertices on~$C$.
Let $z_{v_i}$ be the row in $C_{\mathcal{Q}}(x)$
corresponding to the vertex~$v_i$.
Then we replace the row
corresponding to~$v_1$ by
$$Z:=z_{v_1} - {w(p_1)}z_{v_2} + {w(p_1)w(p_2)}z_{v_3}
-+\ldots \hskip4cm$$
$$\hskip3cm \ldots+ (-1)^{k-1} {w(p_1)\cdots w(p_{k-1})} z_{v_{k}}.
$$
Again, as in the proof of Proposition \ref{prop-finite} all
contributions coming from
paths starting in~$v_1$ with~$p_1$ cancel in~$Z$, and the
contributions from the other paths
starting in~$v_1$ (now there may be nontrivial ones)
occur with the factor
$$1+(-1)^{k-1} {w(p_1)\cdots w(p_{k-1})} \cdot {w(p_k)} =
1-(-1)^{l(C)} {w(C)} \:.$$
At this point our proof has to deviate from the previous one
as there might well be
non-zero paths from vertices other than $v_1$ involving~$p_1$.
(This is because in our
situation we can not guarantee that there is only one incoming
arrow into~$v_1$. Note that in the previous proof we could only
assume this because there were no cycles with no relations.
In the present proposition the quiver might, for instance,
have the property that all vertices are of valency~4.)

In the next step we try to get rid of this problem by performing
column operations. By assumption there exists an
incoming arrow~$q_1$ not on the cycle~$C$.
Consider the unique maximal path~$q$ going backwards from~$q_1$
along the zero relations, i.e. $q=q_l\ldots q_2q_1$ with
$q_jq_{j-1}\in I$ for all $j=l,\ldots,2$. Note that $q_1$
might also belong to an oriented cycle with full relations, namely
if there are two such cycles attached to~$v_1$.
Then $q$ is set to be one walk around the cycle.
Let $v_l',\ldots,v_2',v_1$ be the vertices on this path~$q$.
Moreover, let $\bar{s}_{v_l'},\ldots,\bar{s}_{v_2'},\bar{s}_{v_1}$ be the
corresponding columns of the modified Cartan matrix $\bar{C}(x)$
(obtained from $C_{{\mathcal Q}}^{w}(x)$ by replacing $z_{v_1}$ with~$Z$).
Then we set
$$S:=s_{v_1} - {w(q_1)}s_{v_2'} + {w(q_1)w(q_2)}s_{v_3'}
-+\ldots \hskip4cm$$
$$\hskip3cm \ldots+ (-1)^{l-1} {w(q_1)\cdots w(q_{l-1})} z_{v_l'}.
$$
Let $\widetilde{C}(x)$ be the matrix obtained from
$C_{\mathcal{Q}}(x)$
by first replacing the row $z_{v_1}$ by~$Z$, and then the
column $s_{v_1}$ by~$S$.

Completely analogous to the argument for the row operation,
the contributions in
$s_{v_1}$ from the paths ending with~$q_1$ are cancelled in~$S$.

If $q_1$ does belong to an oriented cycle $C'$ with full relations,
then the contributions from the other paths ending in~$v_1$
occur in~$S$ with a factor
$$1+(-1)^{l-1} {w(q_1)\cdots w(q_{l-1})} \cdot {w(q_l)} =
1-(-1)^{l(C')} {w(C')} \:.$$

If $q_1$ does not belong to an oriented cycle with full relations,
then this factor does not occur.

Let $\bar{\bar{C}}(x)$ be the matrix obtained from $\widetilde{C}(x)$
by taking out the factor $1-(-1)^{l(C)} {w(C)} $ from the first
row
and, if $q_1$ belongs to an oriented cycle with full relations,
the factor $1-(-1)^{l(C')} {w(C')} $ from the first column.

Recall that in $\bar{\bar{C}}(x)$ all contributions from paths starting
with~$p_1$ and all contributions from paths ending with~$q_1$
are cancelled.
But, in general, there will be paths between vertices
other than $v_1$ involving the (non-zero) product~$q_1p_1$.
Therefore, we have to introduce the new arrow~$\bar{p}$ replacing~$q_1p_1$.

Then $\bar{\bar{C}}(x)$ is precisely the weighted Cartan matrix
$C_{\bar{\mathcal{Q}}}^{\bar w}(x)$ of the $\bar w$-weighted
locally gentle quiver described in the proposition.

Summarizing the above arguments we get for the determinants
(leave out the factor $1-(-1)^{l(C')} {w(C')}$ if $q_1$ does not
belong to an oriented cycle with full relations)
\begin{eqnarray*}
\det C_{\mathcal{Q}}^w(x)
 & = & (1-(-1)^{l(C)}{w(C)})\cdot \det\bar{C}(x) \\
 & = & (1-(-1)^{l(C)}{w(C)} )(1-(-1)^{l(C')}{w(C')} )\cdot
\det\bar{\bar{C}}(x) \\
 & = & (1-(-1)^{l(C)}{w(C)} )(1- (-1)^{l(C')}{w(C')} )\cdot
\det C_{\bar{\mathcal{Q}}}^{\bar w}(x),
\end{eqnarray*}
as claimed.
\qed

\begin{remark} \label{rem-reduce}
{\em The proof of Proposition~\ref{prop-reduce} also works
when there is no incoming arrow $q_1$;
in this situation, the quiver $\bar{\mathcal Q}$ is obtained
from ${\mathcal Q}$ by just removing~$p_1$.}
\end{remark}

\subsection{An explicit example} \label{explicit-example}
Let $\mathcal{Q}=(Q,I)$ be
the following
locally gentle quiver
(where relations are indicated by dotted lines).


\begin{center}
\includegraphics[scale=0.75]{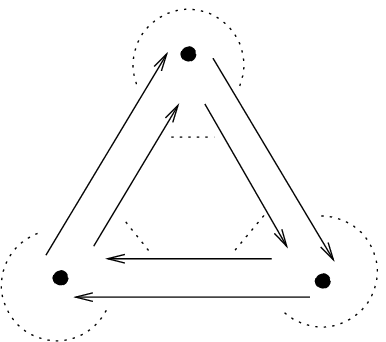}
\end{center}


Note that we have two minimal cycles of length 3 with full relations,
and one minimal cycle with no relations of length~6.
As it is somewhat cumbersome to write out the weighted
Cartan matrix for the generic weight function,
we choose the weight function $w$ to be $q$ on every arrow.
Then the  weighted Cartan matrix (i.e., in this case the $q$-Cartan
matrix) has the form
$$C^w_{\mathcal{Q}}(q) =C_{\mathcal{Q}}(q) = \left(
\displaystyle
\begin{array}{ccc}
\d \frac{1+q^3}{1-q^3} & \d \frac{2q}{1-q^3} &
\d \frac{2q^2}{1-q^3} \\[8pt]
\d \frac{2q^2}{1-q^3} & \d\frac{1+q^3}{1-q^3} &
\d\frac{2q}{1-q^3} \\[8pt]
\d\frac{2q}{1-q^3} & \d\frac{2q^2}{1-q^3} &
\d\frac{1+q^3}{1-q^3}
\end{array}
\right)\:.
$$
As an illustration we go through the corresponding reduction steps
described in~\ref{reduce-ZC} in detail.
We denote the three outer arrows in the above quiver by~$a$, and
the three inner arrows by~$b$.

First, we perform row operations along the cycle~$a^3$ with full
relations, replacing the first row
by $Z:=z_1-qz_2+q^2z_3$. Then we perform column operations
backwards along the cycle~$b^3$,
replacing the first column by the linear combination
$S:=s_1-qs_3+q^2s_2$.
We get the new matrix (having the
same determinant) of the form
$$\widetilde{C}(q) = \left(
\begin{array}{ccc}
\d\frac{1+q^3}{1-q^3} & \d\frac{qt(1+q^3)}{1-q^3} &
\d\frac{q^2(1+q^3)}{1-q^3} \\[8pt]
\d\frac{q^2(1+q^3)}{1-q^3} & \d\frac{1+q^3}{1-q^3} &
\d\frac{2q}{1-q^3} \\[8pt]
\d\frac{q(1+q^3)}{1-q^3} & \d\frac{2q^2}{1-q^3} &
\d\frac{1+q^3}{1-q^3}
\end{array}
\right)
$$
When computing the determinant, we
can now extract the factor $1+q^3$ from the first row and from
the first column and get
$$\det C^w_{\mathcal Q}(q) = (1+q^3)^2\cdot
\det \left(
\begin{array}{ccc}
\d\frac{1}{1-q^6} & \d\frac{q}{1-q^3} &
\d\frac{q^2}{1-q^3} \\[8pt]
\d\frac{q^2}{1-q^3} & \d\frac{1+q^3}{1-q^3} &
\d\frac{2q}{1-q^3} \\[8pt]
\d\frac{q}{1-q^3} & \d\frac{2q^2}{1-q^3} &
\d\frac{1+q^3}{1-q^3}
\end{array}
\right)
$$

The next step in the proof of Proposition~\ref{prop-reduce}
is the transition to the modified
$\widetilde{w}$-weighted quiver
$\widetilde{\mathcal{Q}} :=(\widetilde{Q},\widetilde{I})$ which
now has one arrow of weight~$q^2$ and takes the form


\begin{center}
\includegraphics[scale=0.5]{loc-weight-3.eps}
\end{center}


with conventions on the weights as in Example~\ref{ex-duality}.

The weighted Cartan matrix of this  locally gentle quiver
with respect to the weights $q^2t$ on the marked arrow and $qt$
on all others
has been considered in Example~\ref{ex-duality};
we call this weight function also~$\tilde w$ and keep in mind
that we will have to specialize $t$ to~1 to obtain the matrix appearing on the
right side above.

The transformed quiver
has no more oriented cycles with full relations. Hence
its dual weighted locally gentle quiver
$\widetilde{\mathcal{Q}}^{\#}$
has no cycles with no relations. Thus we can compute
its weighted Cartan determinant from~\ref{finite}:
$$\det C^{\tilde w}_{\widetilde{\mathcal{Q}}^{\#}}(q,t) = 1+q^6t^5.$$
From our duality result Proposition~\ref{duality-thm}
we deduce that
$$\det C^{\tilde w}_{\widetilde{Q}}(q,1) =
(\det C^{\tilde w}_{\widetilde{Q}^{\#}}(q,-1))^{-1} = \frac{1}{1-q^6}.$$
Summarizing the above steps we get
\begin{eqnarray*}
\det C^w_{\mathcal{Q}}(q) & = & \det \widetilde{C}(q)~~ = ~~
(1+q^3)^2\cdot \det C^{\tilde w}_{\widetilde{\mathcal{Q}}}(q,1) \\
 & = & \frac{(1+q^3)^2}{1-q^6}.
\end{eqnarray*}
Note that this is exactly in line with our main theorem
since ${\mathcal Q}=(Q,I)$ has two minimal oriented cycles with full
relations
of length~3 and
one minimal oriented cycle with no relations of length~6.

\subsection{Almost finishing the proof using duality} \label{dual}

We are now going to complete the proof of our main theorem
in most cases. We will encounter some very special locally
weighted gentle quivers which have to be treated separately
in the next subsection.

Let $\mathcal{Q}=(Q,I)$ be an arbitrary locally gentle
quiver with the generic  weight function~$w$.
Let $\mathcal{ZC}(\mathcal{Q})$ be its set of minimal
oriented cycles
with full relations, and let $\mathcal{IC}(\mathcal{Q})$
be its set of minimal cycles with no relations

Assume $\mathcal{ZC}(\mathcal{Q})\neq\emptyset$.
If there is a zero cycle having a vertex that is not incident to four of
the arrows of the cycle, then
we can construct a new
weighted locally gentle quiver $\bar{\mathcal{Q}}$ which
has $|\mathcal{ZC}(\bar{\mathcal{Q}})|<
|\mathcal{ZC}(\mathcal{Q})|$
according to the combinatorial rule described in
Proposition~\ref{prop-reduce}.
(Actually the difference is 1 or 2,
depending on whether the chosen vertex $v_1$ is attached to one or
two minimal oriented cycles with full relations.)
The crucial observation is that in this construction
the number of minimal cycles with no relations and their
weights are not changed at all.
In particular, there is a weight-preserving bijection
$\mathcal{IC}(\mathcal{Q}) \to
\mathcal{IC}(\bar{\mathcal{Q}})$.

By Proposition~\ref{prop-reduce} and Remark~\ref{rem-reduce} each oriented cycle~$C$ with
full relations lost in this transition
from $\mathcal{Q}$ to $\bar{\mathcal{Q}}$ gives a factor
$1-(-1)^{l(C)} {w(C)} $ for the computation of the determinant.

We assume now that on each cycle with full relations
we find a vertex for which
we can perform this reduction step (the only critical
situation occurs when there is no such vertex on each
cycle with full relations in the quiver and its dual, and
we can not reduce further; this will be dealt with in the next
subsection).
Then, continuing this combinatorial process inductively, we reach a
$\tilde w$-weighted locally gentle quiver $\widetilde{\mathcal{Q}}$
with $\mathcal{ZC}(\widetilde{\mathcal{Q}})=\emptyset$
and
\begin{eqnarray*}
\det C_{\mathcal{Q}}^w(x) =
\big(\prod_{C\in \mathcal{ZC}(\mathcal{Q})}
(1- (-1)^{l(C)}{w(C)} )\big)\cdot
\det C_{\widetilde{\mathcal{Q}}}^{\tilde w}(x).
\end{eqnarray*}
Since $\mathcal{ZC}(\widetilde{\mathcal{Q}})=\emptyset$ we will
have for the dual weighted locally gentle quiver that
$\mathcal{IC}(\widetilde{\mathcal{Q}}^{\#})=\emptyset$.
Thus
applying Proposition~\ref{prop-finite}
to $\widetilde{\mathcal{Q}}^{\#}$
(specialize the generic weight function to  $\tilde w$),
duality (via Remark~\ref{rem-dual})
and the fact that there are weight-preserving bijective
correspondences
$$\mathcal{ZC}(\widetilde{\mathcal{Q}}^{\#}) \leftrightarrow
\mathcal{IC}(\widetilde{\mathcal{Q}}) \leftrightarrow
 \mathcal{IC}(\mathcal{Q})
$$
we obtain
\begin{eqnarray*}
\det C^{\tilde w}_{\widetilde{\mathcal{Q}}}(x)
& = & \left( \prod_{C\in\mathcal{IC}(\mathcal{Q})}
(1-{w(C)})\right)^{-1}.
\end{eqnarray*}
Finally, combining the above equations we can conclude that
\begin{eqnarray*}
\det C^w_{\mathcal{Q}}(x)
 & = & \frac{\prod_{C\in \mathcal{ZC}(\mathcal{Q})}
(1- (-1)^{l(C)}{w(C)} )}{\prod_{C\in\mathcal{IC}(\mathcal{Q})}
(1-{w(C)})}.
\end{eqnarray*}
\qed

\subsection{Critical quivers} \label{critical}

The previous reduction steps yield a proof of our main theorem,
unless the $w$-weighted locally gentle quiver ${\mathcal Q}=(Q,I)$
is connected and has the following
property:
each vertex has valency~4, and all the arrows incident to it belong to
the same cycle with no relations and also to the same cycle
with full relations.
Note that then
$Q$ has exactly one minimal cycle with no relations and
exactly one minimal cycle with full relations, and
every arrow of $Q$ belongs to both these cycles.
In other words, the quiver consists of these two cycles, which
are interwoven and both 'eight-shaped' at each vertex.
We call these weighted locally
gentle quivers {\em critical}.
Here are two examples of critical locally gentle quivers.


\begin{center}
\includegraphics[scale=0.55]{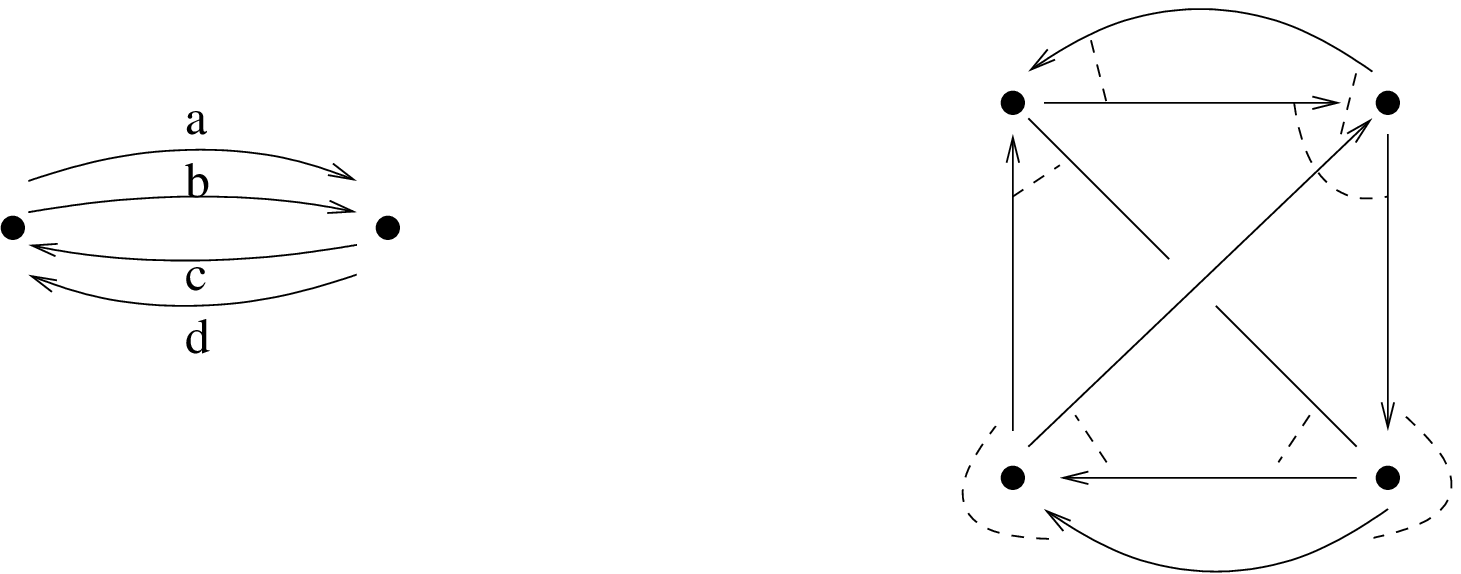}
\end{center}

where for the left quiver we have the relations $ad=0$, $bc=0$,
$db=0$ and $ca=0$. In the right quiver, the zero relations are
indicated by the dotted lines.
\smallskip

In the next subsection below we will discuss a combinatorial interpretation
of critical locally weighted quivers. A consequence of a deep combinatorial
result, the Harer-Zagier formula~\cite{Harer-Zagier}, then implies that such
quivers exist only with an even number of vertices.

Note that the class of
critical locally gentle quivers is closed under duality (as
introduced in Section~\ref{sec-duality}), and that they don't
satisfy the assumptions of the main reduction step
Proposition~\ref{prop-reduce}.
So we indeed have to deal with these critical quivers
separately.

The following result then completes the proof of our main theorem.

\begin{prop} \label{prop-critical}
Let ${\mathcal Q}=(Q,I)$ be a critical
locally gentle quiver
with the generic weight function~$w$.
Then for the weighted Cartan matrix
the following holds
$$\det C^w_{\mathcal Q}(x) = 1.$$
\end{prop}

Note that this is precisely the value of the determinant predicted by
the main theorem. In fact, a critical quiver ${\mathcal Q}$ has precisely
one cycle in $\mathcal{ZC(Q)}$ and one cycle in $\mathcal{IC(Q)}$, both
of which have length $2|Q_0|$ and the same weight (namely the product
of the weights over all arrows in the quiver).
\smallskip

For the proof of the above result we shall need the following very simple
fact about determinants of matrices which only differ in one entry.

\begin{lemma} \label{lemma-det}
Let $C=(c_{ij})$ and $\tilde{C}=(\tilde{c}_{ij})$ be
$n\times n$-matrices. Assume $\tilde{c}_{ij}=c_{ij}$ except for
$i=j=1$. Then
$$\det \tilde{C} = \det C + (\tilde{c}_{11}-c_{11})\det {\mathcal C}_{11}$$
where ${\mathcal C}_{11}$ is the principal submatrix of
$C$ (and of $\tilde{C}$) obtained by removing the first row and
column. \qed
\end{lemma}

\begin{proof} (of Proposition \ref{prop-critical})
Let ${\mathcal Q}=(Q,I)$ be a critical quiver as defined above,
with generic weight function $w$.
We need to introduce some notation. Let $n$ denote the number
of vertices of~$Q$.
Fix any vertex $v_1$
in $Q$, and let $p_1$ be one of the arrows starting in~$v_1$.
Since ${\mathcal Q}$ is critical, there is a unique minimal
cycle $p$ with full relations, starting with~$p_1$,
and containing each arrow of~$Q$
precisely once, and passing through each vertex of~$Q$ twice.
Note that the path~$p$
has length~$2n$.
We can write $p=\bar{p}p'$ where $\bar{p}$ denotes the
initial {\em proper} subpath of~$p$ of positive length
ending in~$v_1$. Let $a$ be the length of~$\bar{p}$.

On the weighted Cartan matrix $C_{{\mathcal Q}}^{w}(x)$
we shall perform row and column operations similar to the
previous proofs. These will again be given by suitable alternating
sums along the cycle $p$. But as we will see, one has to be careful
since these will not be elementary row and column operations, since
each vertex occurs twice on $p$.

More precisely, let $v_1,v_2,\ldots,v_{a},v_{a+1}=v_1,v_{a+2},\ldots,
v_{2n-1},v_{2n}=v_1$ be the vertices, and
$p_1,p_2,\ldots,p_{2n}$ the arrows on $p$. In particular, we have
$\bar{p}=p_1\cdots p_a$.

The row of the Cartan matrix
corresponding to the vertex $v_i$ is denoted by~$z_{v_i}$.
Then consider the following linear combination of
rows
$$Z:=z_{v_1} - {w(p_1)}z_{v_2} + {w(p_1p_2)}z_{v_3} -+ \ldots
+ (-1)^{a-1} {w(p_1\ldots p_{a-1})}z_{v_{a}} +
$$
$$
+ (-1)^{a} {w(p_1\ldots p_{a})}z_{v_{a+1}} +\ldots +
(-1)^{2n-1} {w(p_1\ldots p_{2n-1})}z_{v_{2n}}.
$$
We now consider the matrix obtained from
$C_{{\mathcal Q}}^{w}(x)$ by replacing
the row corresponding to $v_1$ by~$Z$. W.l.o.g.\ we assume
that $z_{v_1}$ is the first row.
Since $v_{a+1}=v_1$, the row $z_{v_1}$ occurs twice in $Z$, namely
with factor 1 and with factor
$(-1)^{a} {w(p_1\ldots p_{a})} = (-1)^{a} {w(\bar{p})}
$. Hence, for the determinant we
obtain

$$\begin{array}{rcl}
\mbox{{\large $\det C_{\mathcal Q}^w(x)$}} & = &
\mbox{{\Large $\frac{1}{1+(-1)^a {w(\bar{p})}}$}}\cdot
\det \mbox{{\footnotesize $\left[
\begin{array}{cccc}
(1-{w(p)})\bar{c}_{11} & \ldots & \ldots & (1-{w(p)})\bar{c}_{1n} \\
c_{21} & \ldots & \ldots & c_{2n} \\
\vdots & & & \vdots \\
c_{n1} & \ldots & \ldots & c_{nn}
\end{array} \right]$}} \\
\end{array}
$$
where $\bar{c}_{1j}$ is the contribution of paths starting
at~$v_1$ but not with the arrow~$p_1$, and ending in~$v_j$.
(Note that indeed, in the alternating sum~$Z$, the contributions
coming from paths starting with~$p_1$ cancel. This is completely
analogous to previous proofs.)

Now we perform column operations. For a vertex~$v_i$, let $s_{v_i}$
denote the column of the above matrix corresponding to~$v_i$.

We consider alternating sums given by
going {\em backwards} along the cycle $p$ with full relations,
starting with~$p_a$. Thus we set
$$S:=s_{v_1} - {w(p_a)}s_{v_a} + {w(p_{a-1}p_a)}s_{v_{a-1}}
+ \ldots +
(-1)^{a-1} {w(p_2\ldots p_{a})} s_{v_2}
+(-1)^{a} {w(\bar{p})}s_{v_{1}}
$$
$$
+ (-1)^{a+1} {w(p_{2n}p_1\ldots p_{a})}s_{v_{2n}} +\ldots +
(-1)^{2n-1} {w(p_{a+2}\ldots p_{2n}p_1\ldots p_{a})}s_{v_{a+2}}.
$$
Then the above determinant becomes
$$
\begin{array}{rcl}
\det C^{w}_{{\mathcal Q}}(x) & = &
\displaystyle
\frac{1-{w(p)}}{(1+(-1)^a {w(\bar{p})})^2}
\cdot \det \mbox{{

$\left(
\begin{array}{cccc}
(1+(-1)^a{w(\bar{p})}) & \bar{c}_{12} & \cdots &
\bar{c}_{1n} \\[5pt]
(1-{w(p)})\bar{c}_{21} & & &  \\
\vdots & & {\mathcal C}_{11} &  \\
(1-{w(p)})\bar{c}_{n1} & & &
\end{array} \right)$}} \\[40pt]
 & = &
\displaystyle
\frac{(1-{w(p)})^2}{(1+(-1)^a {w(\bar{p})})^2}
\cdot \det \mbox{{
$\left( \begin{array}{cccc}
\d \frac{1+(-1)^a {w(\bar{p})}}{1-{w(p)}} &
\bar{c}_{12} & \cdots & \bar{c}_{1n} \\[10pt]
\bar{c}_{21} &  &  & \\
\vdots & & {\mathcal C}_{11} &  \\
\bar{c}_{n1} & & &
\end{array} \right) $}}
\end{array}
$$
where $\bar{c}_{j1}$ is the contribution of paths ending
at~$v_1$ but not with the arrow~$p_{a}$, and starting in~$v_j$.
(Note that indeed, in the alternating sum~$S$, the contributions
coming from paths ending with~$p_a$ cancel. This is again completely
analogous to previous proofs.)

The crucial fact to observe now is that the latter matrix,
apart from the top left entry, is
the weighted Cartan matrix of the $\bar w$-weighted locally gentle
quiver~$\bar{{\mathcal Q}}$
obtained from ${\mathcal Q}$ by removing $p_1$ and~$p_a$, and
replacing them by a new arrow $\pi$ from $v_a$ to~$v_2$ of weight
$\bar w(\pi)=w(p_1)w(p_a)$; all other weights are kept the same
for the weight function~$\bar w$.

By Lemma~\ref{lemma-det} we thus get
$$\begin{array}{rcl}
\det C^{w}_{{\mathcal Q}}(x) & = &
\displaystyle
\frac{(1-{w(p)})^2}{(1+(-1)^a {w(\bar{p})})^2}
\; \cdot \\[8pt]
&&
\displaystyle
\quad \cdot \left( \det C^{\bar w}_{\bar{{\mathcal Q}}}(x) +
\left(\frac{1+(-1)^a{w(\bar{p})}}{1-{w(p)}} -
\frac{1}{1-{w(p)}}\right)
\det {\mathcal C}_{11} \right). \\
\end{array}
$$
The quiver $\bar{{\mathcal Q}}$ has the property that
the vertex $v_1$ only has valency~2. In particular,
$\bar{{\mathcal Q}}$ is not critical, and we can
apply our previous reduction steps which prove
the main theorem for non-critical quivers.
Note that $\bar{{\mathcal Q}}$ has the same minimal cycle without relations
as~${\mathcal Q}$. But the cycle with full relations is broken up,
and the only minimal cycle with full relations in $\bar{{\mathcal Q}}$
is $p_2\ldots p_{a-1}q$, of length $a-1$
and its $\bar w$-weight is equal to the
weight~$w(\bar{p})$.
Thus by our main theorem for non-critical quivers we get
$$\det C^{\bar w}_{\bar{{\mathcal Q}}}(x) =
\frac{1-(-1)^{a-1} {w(\bar{p})}}{1-{w(p)}}.$$

Moreover, the principal minor $\det {\mathcal C}_{11}$ is the weighted
Cartan matrix of the quiver $\bar{{\mathcal Q}}'$ obtained from
${\mathcal Q}$ by removing the vertex $v_1$ (and all arrows
attached to it), and then replacing the non-zero product
$p_ap_1$ by a new arrow of weight $w(p_a)w(p_1)$,
and the non-zero product $p_{2n}p_{a+1}$ by a new arrow
of weight $w(p_{2n})w(p_{a+1})$; let us call the corresponding new
weight function~$w'$. Note that the
quiver $\bar{{\mathcal Q}}'$ has the same cycle with no relations
as~${\mathcal Q}$, but it now has two minimal cycles with full
relations, namely one of length $a-1$ and $w'$-weight $w(\bar{p})$, and
the other of length $2n-a-1$ and $w'$-weight~$w(p')$.
By induction on the number of vertices, we get for the
Cartan determinant of
the $w'$-weighted locally gentle quiver $\bar{{\mathcal Q}}'$ that
$$\det {\mathcal C}_{11} = \det C^{w'}_{\bar{{\mathcal Q}}'}(x) =
 \frac{(1-(-1)^{a-1}{w(\bar{p})})(1-(-1)^{2n-a-1}{w(p')}
)}{1-{w(p)}}.
$$
We can now plug in this information into the above equations for
the Cartan determinant of ${\mathcal Q}$ to get
$$\begin{array}{rcl}
\det C_{{\mathcal Q}}^{w}(x) & = &
\displaystyle
\frac{(1-{w(p)})^2}{(1+(-1)^a{w(\bar{p})})^2}
\left(
\frac{1-(-1)^{a-1}{w(\bar{p})}}{1-{w(p)}} +
\frac{(-1)^a{w(\bar{p})}}{1-{w(p)}}\cdot
\det {\mathcal C}_{11}
\right) \\

& & \\

& = &
\displaystyle
\frac{1-{w(p)}}{1+(-1)^a{w(\bar{p})}}
\left( 
1+
\frac{(-1)^a{w(\bar{p})}(1+(-1)^{2n-a}{w(p')})}{1-{w(p)}}
\right) \\

& & \\

& = &
\displaystyle
\frac{1-{w(p)}}{1+(-1)^a{w(\bar{p})}}
\left(
\frac{1-{w(p)}+(-1)^a{w(\bar{p})}
+(-1)^{2n}{w(p)}}{1-{w(p)}}
\right) \\

& & \\

& = & 1\,,
\end{array}
$$
as claimed.
\end{proof}

\subsection{Critical quivers and combinatorial configurations}
\label{configurations}

Let $\mathcal{Q}$ be a critical locally gentle  quiver.
Then $\mathcal{Q}$ is connected,
and all arrows belong to a single oriented cycle of length~$2n$, and they also
all belong to a single oriented cycle with full relations
of length~$2n$, where $n=|Q_0|$.
The quiver may then also be described in a different way as follows.
We label the vertices from $1$ to $n$ and start walking
on the oriented cycle at vertex~$1$, not repeating any arrow;
this gives a (circular) sequence of length~$2n$, where each
number $1, \ldots,n$ appears twice.
Note that no two consecutive numbers
on this circular sequence are equal, since otherwise
there would be a loop at the corresponding  vertex,
a situation excluded by our condition on the quiver.
Thus we may visualize this by an oriented $2n$-polygon
with $n$ secants, each connecting two vertices with the same label.
The walk along a path with zero relations in the quiver
corresponds in this picture to a walk of the following type: take a
step on the polygon, then slide along the secant to the vertex with
the same label, take again a step on the polygon, then go over the secant
and so on; let us call this a {\em secant walk}.
Indeed, the secants correspond exactly to the dotted lines
indicating the zero relations in our quiver pictures.
Such secant configurations in $2n$-polygons have appeared
also in other contexts, sometimes
in a slightly disguised form,
e.g., see   \cite{Goulden-Nica}, \cite{Harer-Zagier}, \cite{Singmaster}.
\\
In fact, we have a more special situation above.
In an arbitrary configuration as above, there will be several
cyclic secant walks which do not cover all arrows of the polygon.
Our critical quivers correspond to configurations where we
have a cyclic secant walk covering all arrows (and secants);
we call these {\em closed} configurations.
The number of labelled configurations of this type can be determined
by using a formula due to Harer and Zagier~\cite{Harer-Zagier},
for which a combinatorial proof was given by Goulden and
Nica~\cite{Goulden-Nica}.
We now explain the connection between our configurations and
the situation in \cite{Goulden-Nica}; first we introduce the notation
from \cite{Goulden-Nica} and state the formula.

In the symmetric group~$S_{2n}$, let $P_n$ denote the conjugacy class of
involutions without fixed points.
We denote by $\gamma=(1\:2\: \ldots 2n-1\: 2n)$ the cyclic shift permutation
in~$S_{2n}$.
Then set $A_n=\{ \mu\gamma\mid \mu \in P_n\}$ (here is a slight
change in comparison with
\cite{Goulden-Nica} in that we take the products with $\gamma$ instead of
with $\gamma^{-1}$, but this does not affect the following counts).
Now let $a_{n,k}$ be the number of permutations in $A_n$ with exactly $k$
cycles in the disjoint cycle representation.
Recall that for any $n$ one sets
$(2n-1)!!:= 1\cdot 3 \cdot \ldots \cdot (2n-3)\cdot (2n-1)$.
With these notations
the formula obtained by Harer and Zagier reads as follows
(see \cite{Goulden-Nica}):

\begin{theorem}\label{HZ} \cite{Harer-Zagier}
For $n\geq 1$,
$$\sum_{k\geq 1} a_{n,k}x^k=(2n-1)!! \sum_{k\geq 1} 2^{k-1}\binom{n}{k-1}
\binom{x}{k}\:.$$
\end{theorem}

From this formula, we may easily derive an explicit formula for $a_{n,1}$:

\begin{cor}
For $n\geq 1$,
$$a_{n,1}=\left\{
\begin{array}{cl}
\displaystyle \frac{(2n-1)!!}{n+1}
& \mbox{if \, $n$\, is even}\\[8pt]
0 &  \mbox{if\, $n$\, is odd}
\end{array}\right.
\:.$$
\end{cor}

Clearly, a secant configuration on a labelled oriented $2n$-polygon
may equivalently be described by an involution by walking along the
vertices $v_1,\ldots,v_{2n}$
on the oriented cycle, and then defining the involution as the product
of all transpositions
$(i\, j)$ with $v_i=v_j$ (i.e., $v_i$ and $v_j$ are joined by a secant).
This is an involution without fixed points, and hence an element in~$P_n$.
But note, that in our secant configurations we never join two
neighbouring vertices,
so we only get $\sigma\in P_n$ with
$\sigma(i)\neq i\pm 1$ (modulo $2n$) for all~$i$,
and indeed, we obtain {\em all} those involutions;
we denote this subset of $P_n$ by~$P'_n$.
\\
Computing the product $\mu\gamma$ for  $\mu \in P'_n$ corresponds
exactly to taking
secant walks in the secant configuration,
i.e., the cycles in this product correspond to the cyclic
secant walks in the $2n$-polygon.
In particular, the configuration corresponding to $\mu$ is closed exactly if
$\mu\gamma$ is a  $2n$-cycle in~$S_{2n}$.
Now note that if $\sigma\in P_n$ with $\sigma(i)= i+1$ for some $i$, then
$i$ is a fixed point of $\sigma\gamma$; hence $\sigma$ does not give a
contribution
to~$a_{n,1}$. This shows that
$a_{n,1}$ is the number of permutations in
$A'_n=\{ \mu\gamma\mid \mu \in P'_n\}$ which are $2n$-cycles.
By the previous discussion, we have thus shown that $a_{n,1}$ is also
the number of closed secant configurations on the labelled oriented
$2n$-polygon,
and hence this is the number of critical quivers at the beginning of
this section.
\medskip

Motivated by the quiver situation, we are even more interested in counting
unlabelled configurations as above (or equivalently, counting the
configurations
on a regular $2n$-polygon up to dihedral symmetry).
Even without the restriction on counting only closed secant configurations
this is a difficult problem, see \cite{Singmaster} where the  values
up to $n=8$
were computed; with somewhat improved methods and today's computers one
can easily
extend this list but a closed formula still does not seem to be known.


\end{document}